\documentclass[a4paper,11pt]{article}

\usepackage{amsmath,amssymb,amsthm}
\theoremstyle{plain}
\newtheorem{thm}{Theorem}
\newtheorem{lm}{Lemma}

\newtheorem{cor}{Corollary}
\usepackage{caption,graphicx}
\usepackage[all]{xy}

  \renewenvironment{thebibliography}[1]{%
    \begin{oldthebibliography}{#1}%
      \setlength{\parskip}{0ex}%
      \setlength{\itemsep}{0ex}%
  }%
  {%
    \end{oldthebibliography}%
  }

\begin{document}












\noindent{Topology/ Algebraic Geometry}

\medskip
\noindent{\large \textbf{Braids, conformal module and entropy}}

\medskip
\noindent{\large \textbf{Tresses, module conforme et entropie}}

\medskip

\noindent{Burglind J\"oricke}

\medskip

\noindent{\small Centre de Recerca Matem\`{a}tica, 
E-08193 Bellaterra (Barcelona), Tel +34 93 581 1081, Fax +34 93 581
2202}

\noindent{\small E-mail: joericke@googlemail.com}
\bigskip

\noindent \textbf{Abstract}. The conformal module of conjugacy
classes of braids implicitly appeared in a paper of Lin and Gorin in
connection with their interest in the 13. Hilbert Problem. This
invariant is the supremum of conformal modules (in the sense of
Ahlfors) of certain annuli related to the conjugacy class. This note
states
that 
the conformal module is inverse
proportional to a popular dynamical braid invariant, the entropy.
The entropy appeared in connection with Thurston's theory of surface
homeomorphisms. An application of the concept of conformal module to
algebraic geometry is given.

\medskip
\noindent  \textbf{R\'{e}sum\'{e}.} Le module conforme des classes de
conjugaison de tresses est apparu implicitement dans un article de
Lin et Gorin dans le cadre de leur int\'{e}r\^{e}t pour le 13e probl\`{e}me de
Hilbert. Cet invariant est la borne sup\'{e}rieure des modules conformes
(dans le sens d'Ahlfors) de certains anneaux associ\'{e}s \`{a} la classe de
conjugaison. Cette note affirme que le module conforme est
inversement proportionnel \`{a} un invariant dynamique bien connu des
tresses, l'entropie. L'entropie est apparue dans le cadre de la
th\'{e}orie de Thurston des hom\'{e}omorphismes de surfaces. Une application
du concept de module conforme \`{a} la g\'{e}om\'{e}trie alg\'{e}brique est donn\'{e}e.

\vspace{0.5cm}

\noindent Braids occur in several mathematical fields, sometimes
unexpectedly, and one can think of them in several different ways.
Braids on $n$ strands can be interpreted as algebraic objects,
namely, as elements of the Artin group $\mathcal{B}_n$, or as
isotopy classes of geometric braids, or as elements of the mapping
class group of the n-punctured disc (\cite{Bi}). In connection with
his interest in the Thirteen's Hilbert Problem Arnol'd gave the
following interpretation of the braid group $\mathcal{B}_n$. Denote
by $\mathfrak{P}_n$ the space of monic polynomials of degree $n$
without multiple zeros. This space can be parametrized either by the
coefficients or by the unordered tuple of zeros of polynomials. This
makes  $\mathfrak{P}_n$ a complex manifold, in fact, the complement
of the algebraic hypersurface $\{D_n=0\}$ in complex Euclidean space
$\mathbb{C}^n$. Here $D_n(p)$ denotes the discriminant of the
polynomial $p$. The function $D_n$ is a polynomial in the
coefficients of $p$. Arnol'd studied the topological invariants of
$\mathfrak{P}_n$ (\cite{A}). Choose a base point $p \in
\mathfrak{P}_n$. Using the second parametrization Arnol'd
interpreted the group $\mathcal{B}_n$ of $n$-braids as elements of
the fundamental group $\pi_1(\mathfrak{P}_n,p)$ with base point $p$.

The conjugacy classes $\widehat{ \mathcal{B}}_n$ of the braid group,
equivalently, of the fundamental group $\pi_1(\mathfrak{P}_n,p)$,
can be interpreted as free isotopy classes of loops in
$\mathfrak{P}_n$. We define a collection of conformal invariants of
the complex manifold $\mathfrak{P}_n$. Consider an element $\hat b$
of $\widehat {\mathcal{B}}_n$. We say that a continuous mapping $f$
of an annulus $A= \{z \in \mathbb{C}: \, r<|z|<R\},\;$ $0\leq r < R
\leq \infty,\;$ into $\mathfrak{P}_n$ represents $\hat b$ if for
some (and hence for any) circle $\,\{|z|=\rho \} \subset A\,$ the
loop $\,f:\{|z|=\rho \} \rightarrow \mathfrak{P}_n\,$ represents
$\,\hat b$.
Ahlfors defined the conformal module of an annulus $\,A= \{z \in
\mathbb{C}:\; r<|z|<R\}\,,\,$ as $\,m(A)=
\frac{1}{2\pi}\,\log(\frac{R}{r})\,.\, $ Two annuli of finite
conformal module are conformally equivalent iff they have equal
conformal module. If a manifold $\Omega$ is conformally equivalent
to an annulus $A$ its conformal module is defined to be $m(A)$.
Associate to each conjugacy class of the fundamental group of
$\mathfrak{P}_n$, or, equivalently, to each
conjugacy class of $n$-braids, its conformal module defined as follows.\\

\noindent  \textbf{Definition} \textit{Let $\hat b$ be a conjugacy
class of $n$-braids, $n \geq 2$. The conformal module $M(\hat b)$ of
$b$ is defined as $ M(\hat b)= sup_{\mathcal{A}}\, m(A),$ where
$\mathcal{A}$ denotes the set of all annuli which admit a
holomorphic mapping into
$\mathfrak{P}_n$ which represents $\hat b$.}\\

\noindent For any complex manifold the conformal module of conjugacy
classes of its fundamental group can be defined. The collection of
conformal modules of all conjugacy classes is a biholomorphic
invariant of the manifold. This concept seems to be especially
useful for locally symmetric spaces, for instance, for the quotient
of the $n$-dimensional round complex ball by a subgroup of its
automorphism group which acts freely and properly discontinuously.
In this case the universal covering is the ball and the fundamental
group of the quotient manifold can be identified with the group of
covering translations. For each covering translation the problem is
to consider the quotient of the ball by the action of the group
generated by this single covering translation and to maximize the
conformal module of annuli which admit holomorphic mappings into
this quotient. The latter concept can be generalized to general
mapping class groups. The generalization has relations to symplectic
fibrations.

Runge's approximation theorem shows that the conformal module is
positive for any conjugacy class of braids. The concept of the
conformal module of conjugacy classes of braids appeared (without
name) in the paper \cite{GL} which was motivated by the interest of
the authors in Hilbert's Thirteen's Problem for algebraic functions.

The following objects related to $\mathfrak{P}_n$ have been
considered in this connection. A continuous mapping of a (usually
open and connected) Riemann surface $X$ into the set of monic
polynomials of degree $n$ (maybe, with multiple zeros) is a
quasipolynomial of degree $n$. It can be written as $\,P(x,\zeta)=
a_0(x) + a_1(x)\zeta + ... + a_{n-1}(x)\zeta^{n-1} + \zeta^n, \, x
\in X,\, \zeta \in \mathbb{C},\; $ for continuous functions $a_j, \;
j=1,...,n\,,$ on $X$. If the mapping is holomorphic it is called an
algebroid function. If the image of the map is contained in
$\mathfrak{P}_n$ it is called separable. A separable quasipolynomial
is called solvable if it can be globally written as a product of
quasipolynomials of degree $1$, and is called irreducuble if it
cannot be written as product of two quasipolynomials of positive
degree. Two separable quasipolynomials are isotopic if there is a
continuous family of separable quasipolynomials joining them. An
algebroid function on the complex line $\mathbb C$ whose
coefficients are polynomials is called an algebraic function. A
quasipolynomial $P$ on $X$ can be considered as a function on
$\,X\times \mathbb{C}\,$. Its zero set $\,\mathfrak{S}_P=
\{(x,\zeta) \in X\times \mathbb{C},\; P(x,\zeta)=0\}\;$ is a
symplectic surface, called braided surface due to its relation to
braids.



The conformal module of conjugacy classes of braids serves as
obstruction for the existence of isotopies of quasipolynomials
(respectively, of braided surfaces) to algebroid functions
(respectively, to complex curves). Indeed,
let $X$ be an open Riemann surface. Suppose $f$ is a quasipolynomial
of degree $n$ on $X$. Consider any domain $A \subset X$ which is
conformally equivalent to an annulus. The restriction of $f$ to $A$
defines a mapping of the domain $A$ into the space of polynomials
$\mathfrak{P}_n$, hence it defines a conjugacy class of
$n$-braids $\hat b_{f,A}$.\\

\begin{lm}
If $f$ is algebroid then $\, mA \leq M(\hat b_{f,A})\,. $
\end{lm}

\noindent Before giving examples of applications of the concept of
the conformal module of conjugacy classes of braids we compare this
concept with a dynamical concept related to braids. Let $\mathbb{D}$
be the unit disc in the complex plane. Denote by $E_n^0$ the set
consisting of the $n$ points $0, \frac{1}{n},...,\frac{n-1}{n}$.
Consider homeomorphisms of the $n$-punctured disc
$\overline{\mathbb{D}} \setminus E_n^0$, which fix the boundary
$\partial{\mathbb{D}}$ pointwise. Equivalently, these are
homeomorphisms of the closed disc $\overline{\mathbb{D}}$ which fix
the boundary pointwise and the set $E_n^0$ setwise. Equip this set
of homeomorphisms with compact open topology. The connected
components of this space form a group, called mapping class group of
the $n$-punctured disc. This group is isomorphic to $\mathcal{B}_n$
(\cite{Bi}). Denote by $\mathcal{H}_b$ the connected component which
corresponds to the braid $b$.

For a homeomorphisms of a compact topological space its topological
entropy is an invariant which measures the complexity of its
behaviour in terms of iterations. It is defined in terms of the
action of the homeomorphism on open covers of the compact space. For
a precise definition of topological entropy we refer to the papers
\cite{AKM} or \cite{F}. For a braid $b$ we define its entropy as
$h(b)=\inf \{h(\varphi): \varphi \in \mathcal{H}_b\}$. The value is
invariant under conjugation with self-homeomorphisms of the closed
disc $\overline{\mathbb{D}}$ 
, hence it does not depend on the position of the set of punctures
and on the choice of the representative of the conjugacy class $\hat
b$. We write $h(\hat b) = h(b)$.

Entropy is a dynamical invariant. It has been considered in
connection with Thurston's theory of surface homeomorphisms.
Thurston himself used dynamical methods (Markov partitions) to show
that many mapping class tori carry a complete hyperbolic metric of
finite volume. Detailed proofs of Thurston's theorems are given in
\cite{F} where also the entropy of homeomorphisms of closed Riemann
surfaces is studied. The study has been extended to Riemann surfaces
with punctures. The common definition for braids is given in the
irreducible case and uses mapping classes of the $n$-punctured
complex plane rather than of the $n$-punctured disc. One can show
that this definition is equivalent to the definition given above.
Entropy has been studied intensively. E.g., the lowest non-vanishing
entropy $h_n$ among irreducible braids on $n$ strands, $n \geq 3,\,$
has been estimated from below by $\frac{\log 2}{4} n^{-1}$
(\cite{P}) and has been computed for small $n$. There is an
algorithm for computing the entropy of irreducible braids
(respectively, of irreducible mapping classes) (\cite{BH}). Fluid
mechanics related to stirring devises uses the entropy of the
arising braids as a measure of complexity.

It turns out that the dynamical aspect and the conformal aspect are
related. The following theorem holds.
\begin{thm} For each $\hat b \in \widehat {\mathcal{B}}_n$ ($n\geq 2)$
$$M(\hat b)=
\frac{\pi}{2}\,\frac{1}{h(\hat b)}\,.$$
\end{thm}
\noindent The proof of the theorem deeply relies on Teichm\"{u}ller
theory, including Royden's theorem on equality of the Teichm\"{u}ller
metric and the Kobayashi metric on Teichm\"{u}ller space, and Bers'
theory of reducible surface homeomorphisms.
\begin{cor} For each $\hat b \in \widehat {\mathcal{B}}_n$ ($n\geq 2)$ and each nonzero
integer $l$
$$M(\widehat{ b\,^l})=\frac{1}{|l|}M(\hat b)\,. $$
\end{cor}


Theorem 1 and \cite{P} allow to re-interpret and to reprove a result
of \cite{Z} on conformal module which we discuss for simplicity only
for prime numbers $n$. Extending results of Lin and Gorin, Zjuzin
\cite{Z} proved that a separable algebroid function $p(z),\, z \in
A,$ of degree a prime number $n$ on an annulus $A$ is reducible if
(a) the conformal module $m(A)$ is bigger than $n \cdot r_0$ for
some absolute constant $r_0$ and (b) the index of the discriminant
$D_n(p(z))$ is divisible by $n$. (The index is the degree of the
mapping $z \rightarrow D_n(p(z))\cdot |D_n(p(z))|^{-1}$ from the
annulus into the circle.) Indeed, if the quasipolynomial is
irreducible then the induced braid class $\hat b_{f,A}$ corresponds
to the conjugacy class of $n$-cycles and, hence, $\hat b_{f,A}$ is
represented by irreducible braids. Condition (a) with $r_0= \frac{2
\pi}{\log 2}$ implies by \cite{P} and theorem 1 that these braids
must be periodic. But condition (b) excludes such periodic braids.

Let $X$ be an open Riemann surface of finite genus with at most
countably many ends. By \cite{HS} $X$ is conformally equivalent to a
domain $\Omega$ on a closed Riemann surface $R$ such that the
connected components of $R \setminus \Omega$ are all points or
closed geometric discs. A geometric disc is a topological disc whose
lift to the universal covering is a round disc (in the standard
metric of the covering). If all connected components of $R \setminus
\Omega$ are points then $X$ is called of first kind, otherwise it is
called of second kind.

The following lemma gives a condition for solvability of algebroid
functions on an annulus.

\begin{lm} Let $X$ be a closed Riemann surface of positive genus with a geometric disc
removed. Suppose $f$ is an irreducible separable algebroid function
of degree $3$ on $X$. Suppose $X$ contains a domain $A$ one of whose
boundary components coincides with the boundary circle of $X$, such
that $A$ is conformally equivalent to an annulus of conformal module
strictly larger than $\frac{\pi}{2}\, (\log( \frac{3
+\sqrt{5}}{2}))^{-1}$. Then $f$ is solvable over $A$.
\end{lm}

\noindent Note that $\log( \frac{3 +\sqrt{5}}{2})$ is the smallest
non-vanishing entropy among irreducible $3$-braids. The estimate of
the conformal module of $A$ and the properties of the covering
$\mathfrak{S}_f \rightarrow X$ allow only the solvable case or the
case of periodic conjugacy classes $\hat b_{f,A}$ which correspond
to conjugacy classes of $3$-cycles. The latter is impossible for
conjugacy classes of products of braid commutators.

Using Lemmas 1 and 2 we give the following application to algebroid
functions on a torus with a hole.

Let $X$ and $Y$ be open Riemann surfaces and let $f$ be a
quasipolynomial on $X$.
Let $w:X\rightarrow Y$ be a homeomorphism. The homeomorphism $w$ can
be interpreted as a new complex structure on $X$. 
Denote by $f_w$ the quasipolynomial
$f_w(y,z)=f(w^{-1}(y),\,\zeta),\; y\in Y,\; \zeta \in \mathbb{C}\,
$. We say that $f$ is isotopic to an algebroid function for the
complex structure $w$ if $f_w$ is isotopic to an algebroid function
on $Y$.

\begin{thm} Let X be a torus with a geometric disc removed. There exist
eight conformal structures of second kind on $X$
with the following property. If $f$ is an irreducible separable
quasipolynomial of degree $3$ on $X$ which is isotopic to an
algebroid function for each of the eight conformal structures on $X$
then $f$ is isotopic to an algebroid function for each complex
structure on $X$ including complex structures of first kind
(determining punctured tori).
\end{thm}

\noindent The conformal structures on the Riemann surface in Theorem
2 are chosen so that each contains a domain conformally equivalent
to an annulus of conformal module estimated from below as in lemma 2
and representing certain element of $\pi_1(X)$.

Let $X$ be as in Theorem 2. Note that the fundamental group
$\pi_1(X,x)$ of $X$ with base point $x \in X$ is a free group on two
generators.
Consider a quasipolynomial of degree $n$ and all its isotopies which
fix the value at the base point $x$. 
For each element $a$ of the fundamental group of $\pi_1(X,x)$ this
defines a homotopy class of loops $\varphi (a)$ with a base point
$p$ in $ \mathfrak{P}_n$, hence an $n$-braid . The mapping $\varphi$
is a homomorphism from $\pi_1(X,x)$ to $\mathcal{B}_n$. There is a
one-to-one correspondence between free isotopy classes of
quasipolynomials (without fixing the value at a point) and conjugacy
classes of homomorphisms from $\pi_1(X)$ to $\mathcal{B}_n$. The
following theorem holds. It shows that the conditions of Theorem 2
are met very rarely.

\begin{thm} Let $X$ be a torus with a hole and $f$ an irreducible separable
quasipolynomial of degree $n$ which is isotopic to an algebroid
function for each complex structure on $X$. Then $f$ corresponds to
the conjugacy class of a homomorphism $\varphi: \pi_1(X) \rightarrow
\Gamma$ where $\Gamma$ is the subgroup of $\mathcal{B}_3$ generated
by the periodic braid $\sigma_1 \sigma_2$. Moreover,
$\varphi(\pi_1(X))$ contains elements other than powers of the
Garside element $(\sigma_1 \sigma_2)^3$.
\end{thm}

\noindent Here  $\sigma_1$ and $ \sigma_2$ are the standard
generators of $\mathcal{B}_3$. Associate to $\sigma_1 \sigma_2$ the
M\"{o}bius transformation $A(\zeta)=e^{\frac{- 2\pi i}{3}} \cdot
\zeta,\; \zeta \in \mathbb{P}^1$. The homomorphism $\varphi$ in
Theorem 3 defines a homomorphism $h$ from $\pi_1(X)$ into the
subgroup of the group of M\"{o}bius transformations which is generated
by $A$. Hence, the fundamental group of any closed torus $T$ acts
freely and properly discontinuously on  $\mathbb{C} \times
\mathbb{P}^1$ by the holomorphic transformations $(z,\zeta) \to
(\gamma (z),h(\gamma)(\zeta)\,),\, $ $\gamma \in \pi_1(X)$. The
quotient is the total space of a holomorphic line bundle over $T$.
The bundle over a punctured torus is trivial and defines a separable
algebroid function $ \tilde f$ in the isotopy class of $f$. The set
$\mathfrak{S}_{\tilde f}$ is a leaf of a holomorphic foliation of
the total space of the bundle over the punctured torus.

There is also the concept of the conformal module of braids rather
than 
of conjugacy classes of braids. This notion is based on the
conformal module of rectangles which admit holomorphic mappings into
$\mathfrak{P}_n$ with suitable boundary conditions on a pair of
opposite sides. Recall that the conformal module of a rectangle with
sides parallel to the coordinate axes is the ratio of the
sidelengths of horizontal and vertical sides. The conformal module
of braids is a finer invariant than entropy. If suitably defined it
is more appropriate for application, in particular, to real
algebraic geometry. In the case of three-braids there are two
versions differing by the boundary conditions on horizontal sides of
rectangles. (Both should be used.) For the first version one
requires that horizontal sides are mapped to polynomials with all
zeros on a real line. In the second version one requires that two of
the zeros have equal distance from the third. These are the cases
appearing on the real axis for polynomials with real coefficients.
The invariant can be studied by quasiconformal mappings and elliptic
functions. The situation for braids on more than $3$ strands is more
subtle. We intend to come back to this concept in a later paper.

\medskip

\noindent \textbf{Acknowledgement}. The author wants to thank
Weizmann Institute, IHES, CRM and the Universities of Grenoble,
Bern, Calais and Toulouse, where the present work has been done. She
is grateful to many mathematicians, especially to F.Dahmani,
P.Eyssidieux, S.Orevkov, O.Viro and M.Zaidenberg for stimulating
discussions.

\bibliographystyle{IEEEtranSA}
\bibliography{}

\end{document}